\documentclass[11pt]{article}
\usepackage{verbatim}
\usepackage{amssymb}
\usepackage{latexsym}
\usepackage{amscd}
\newtheorem{theo}{Theorem}[section]

\newtheorem{lem}[theo]{Lemma}

\newtheorem{quest}[theo]{Question}

\newtheorem{em-bez}[theo]{Notation}

\newtheorem{em-ex}[theo]{Examples and Remarks}

\newtheorem{em-remarks}[theo]{}

\newtheorem{em-re}[theo]{Remark}

\newtheorem{em-defi}[theo]{Definition}

\newenvironment{proof}{\noindent {\sc Proof. }}  {~\hfill$\Box$\par \medskip}
\usepackage{amsmath}

\def\A{{\cal A}}

\def\S{{\cal S}}

\def\NN{{\mathbb N}} 

\def\QQ{{\mathbb Q}}
\def\RR{{\mathbb R}}
\def\TT{{\mathbb T}}
\def\ZZ{{\mathbb Z}}

\def\uu{{\mathbf{u}}}
\def\dd{{\mathbf{d}}}
\def\aa{{\mathbf{a}}}

\title{A note on characterized and statistically characterized subgroups of $\TT=\RR/\ZZ$}
\author{Hans Weber
\\{\footnotesize \tt hans.weber@uniud.it }
\\{\footnotesize  33100 Udine, Italy}}

\begin{document}

\date{}

\maketitle



\begin{abstract}
P. Das, A. Ghosh and T. Aziz has given in \cite[Theorem 3.15]{ADG} a result on statistically characterized subgroups of the circle group $\RR/\ZZ$, which answers, together with \cite[Corollary 2.4]{DG24}, questions of \cite{DDB} and \cite{DG25}. Here we give a completely different and much shorter proof of these results.
\end{abstract} 

\section{Introduction}\label{}
We first fix some notations. Let $\TT:=(\RR,+)/\ZZ$ be the circle group and $\varphi :\RR\rightarrow\TT$ the quotient map. We will use the group-seminorm on $\RR$ defined by $\Vert x\Vert:=\min_{k\in \ZZ}|x-k|$. Then $\Vert x_n\Vert \rightarrow 0$ iff $\varphi(x_n)\rightarrow 0$ in $\TT$ for a sequence $(x_n)$ in $\RR$. 
Let $\S$ be the set of all strictly increasing sequences in $\NN=\{1,2,3,\dots\}$  and $\A$ its subfamily of sequences $(a_n)$ such that $a_{n}|a_{n+1}$ for all $n\in\NN$.   
For $\uu=(u_n)_{n\in\NN}\in\S$ the subgroup 
$$t_{\uu}(\TT):=\{t\in\TT: (u_nt) \text{ converges to 0 in }\TT\}=
\{\varphi(x): x\in\RR\,,\, \Vert u_n x\Vert \text{ converges to 0 in }\RR\}$$ 
of $\TT$ is called, in the terminology of \cite{BDS}, the subgroup of $\TT$ characterized by $\uu$. For the significance of these subgroups of $\TT$ we refer to the survey paper \cite{DDG}, see also \cite{BDMW08}.

As observed in the introduction of \cite{BDMW08}, $t_{\uu}(\TT)$ is either countable or has size $\mathfrak{c}:=2^{\aleph_0}$. In \cite{BDMW03} it is proved:

\begin{theo}\label{In1} \cite[Theorems 3.1 and 3.3]{BDMW03}

Let $\uu=(u_n)\in\S$, $u_0:=0$ and $q_n:=u_{n}/u_{n-1}$ for $n\in\NN$.

(a) If $q_n\rightarrow +\infty$, then $|t_\uu(\TT)|=\mathfrak{c}$.

(b) If $(q_n)$ is bounded, then $t_\uu(\TT)$ is countable.
\end{theo}

If $\uu\in\A$, one can get more information on the structure of $t_{\uu}(\TT)$.

\begin{theo}\label{In2} \cite[Corollary 2.8]{DD} or \cite[Corollary 3.8]{DI} or \cite[Theorem 5.1]{BDMW08}

Let $\uu=(u_n)\in\A$, $u_0:=0$ and $q_n:=u_{n}/u_{n-1}$ for $n\in\NN$.

Then $(q_n)$ is bounded iff $t_\uu(\TT)$ is countable iff $t_\uu(\TT)$ is a torsion group.
\end{theo}

Since for any $\uu\in \S$ the torsion subgroup of $t_\uu(\TT)$ can easily be described (see \cite[Section 2.1]{BDMW08}), together with Theorem \ref{In2} one gets a precise description of $t_\uu(\TT)$ if $\uu\in\A$ and 
$u_{n+1}/u_n$ is bounded:

\begin{theo}\label{bounded}
Let $\aa=(a_n)	\in\A$ such that the sequence $a_{n+1}/a_{n}$ is bounded. 

Then $t_\aa(\TT)=\varphi (\langle\{\frac{1}{a_n}:n\in\NN\}\rangle)$.
\end{theo}

Since for a general sequence $\uu\in\S$ it is difficult to describe $t_\uu(\TT)$, of particular interest, also for examples and counterexamples,  is the sequence $\dd=(d_n)$ defined as follows:

Let $(a_n)\in\A$, $q_n:=\frac{a_n}{a_{n-1}}$ where $a_0:=1$, and let $(d_n)\in\S$ be the sequence defined by 
$$\{d_n:n\in\NN\}=\{ra_k : r,k\in\NN ,  1\leq r<q_{k+1}\}.\quad\quad\quad (\sharp) $$ 

D. Dikranjan and K. Kunen  \cite{DK} showed (see also the revised proof in \cite[Section 3]{B}):

\begin{theo} \cite[Proposition 1.3]{DK} Write $(\zeta_n)$ instead of $(d_n)$ in the special case that $a_n=n!$. Then  $t_{(\zeta_n)}(\TT)=\QQ/\ZZ$.
\end{theo}

In \cite{DG24}, P. Das and A. Ghosh proved the following interesting generalization of this result:

\begin{theo}\label{an-dn}\cite[Corollary 2.4]{DG24}

Let $(a_n)\in\A$ and $\dd=(d_n)\in\S$  defined by  ($\sharp$).
Then $t_\dd(\TT)=\varphi (\langle\{\frac{1}{a_n}:n\in\NN\}\rangle)$.
\end{theo}

In Section \ref{Ch} we provide a very short proof for a theorem with generalizes both, Theorem \ref{bounded} and Theorem \ref{an-dn}.

Recently Dikranjan et al. \cite{DDB} introduced the subgroup $t_{\uu}^s(\TT)$ of $\TT$, \textit{statistically characterized by} $\uu$, replacing in the definition of $t_{\uu}(\TT)$ convergence by statistical convergence, i.e., $t_{\uu}^s(\TT):=\{t\in\TT: (u_nt) \text{ converges statistically to 0 in }\TT\}$. Recall that a sequence $(x_n)$ in a metric space $(X,\rho)$ converges statistically to $x_0$ if $\bar{d}(\{n\in\NN: \rho(x_n,x_0)\geq \varepsilon\}=0$ for every $\varepsilon>0$. Here 
$\bar{d}(A):=\limsup_{n\rightarrow\infty} \frac{|A\cap [1,n]| }{n}$ denotes the upper natural density of $A\subseteq\NN$. 

In \cite{DDB} it is proved:

\begin{theo}\label{BC}\cite[Theorems B and C]{DDB}

If $\aa\in\A$, then $|t^s_\aa(\TT)|=\mathfrak{c}$ and $t^s_\aa(\TT)\neq t_\aa(\TT)$.
\end{theo}

In view of this result the following natural question arises:

\begin{quest}\label{Q1}\cite[Question 6.3]{DDB}

Is it true that $|t^s_\uu(\TT)|=\mathfrak{c}$ and $t^s_\uu(\TT)\neq t_\uu(\TT)$ for any $\uu\in\S$?
\end{quest}

Related to this question is the following more specific question:

\begin{quest}\label{Q3}\cite[2.16]{DG25}

Is $|t^s_\dd(\TT)|=\mathfrak{c}$ for any $\aa\in\A$ where $\dd$ is defined by $(\sharp)$?
\end{quest}

These questions are answered by P. Das, A. Ghosh and T. Aziz \cite{ADG}: In \cite[Theorem 3.15]{ADG} it is given a condition under which $t^s_\dd(\TT)=t_\dd(\TT)$, and in \cite[Corollary 2.4]{DG24} it is proved that $t_\dd(\TT)$ is countable (cf. Theorem \ref{T2} below).

The main aim of this article is to give a completely different and much shorter proof of these results. 

\section{Characterized subgroups of $\TT$}\label{Ch}
The following lemmata are used in the proof of Theorem \ref{T0} as well as in Section \ref{sCh}.
\begin{lem}\label{basic0}
Let $z\in\RR$ and $v\in\NN$. Then $\Vert vz\Vert=\Vert v\Vert z\Vert\Vert$.
If $v\Vert z\Vert\leq\frac{1}{2}$, then $\Vert vz\Vert=v\Vert z\Vert$.
\end{lem}

\begin{proof}
Let $k\in\ZZ$ and $\alpha\in [-\frac{1}{2},\frac{1}{2}]$ with $z=k+\alpha$. Then $\Vert z\Vert=|\alpha|$  and $vz\equiv_\ZZ v\alpha$, therefore $\Vert vz\Vert=\Vert v\alpha\Vert=\Vert v\Vert z\Vert\Vert$.
If $v \Vert z\Vert\leq\frac{1}{2}$, then $\Vert v \Vert z\Vert\Vert=v\Vert z\Vert$.
\end{proof}

\begin{lem}\label{basic1}
Let $(v_n)\in\S$ , $v_0=1$ and $\sup_{n\in\NN}\frac{v_n}{v_{n-1}}\leq q<\infty$. 

Let $z\in\RR$ with $0<\Vert z\Vert\leq\gamma= \frac{1}{2q}$. Then there exists $m\in\NN$ with 
$\Vert v_mz\Vert>\gamma$.
\end{lem}

\begin{proof}
Let $m:=\min\{i\in\NN: v_i\Vert z\Vert >\gamma\}$. Then $v_{m-1}\Vert z\Vert\leq\gamma < v_m\Vert z\Vert$. Therefore $v_m\Vert z\Vert\leq qv_{m-1}\Vert z\Vert\leq q\gamma =\frac{1}{2}$, hence   
$\gamma<v_m\Vert z\Vert=\Vert v_mz\Vert$ by Lemma \ref{basic0}.
\end{proof}

\begin{theo}\label{T0}
Let $\uu=(u_n)\in\S$  and $q_n=\frac{u_n}{u_{n-1}}$ where $u_0:=1$.
Let $a_k=u_{n_k}$ be a subsequence of $(u_n)$ such that $a_k|u_i$ for $i,k\in\NN$ with $i\geq n_k$.

If $(q_n)$ is bounded, then $t_\uu(\TT)=\varphi (\langle\{\frac{1}{a_n}:n\in\NN\}\rangle)$ is the union of the finite subgroups 
$\langle\varphi(\frac{1}{a_n})\rangle$ of $\TT$.
\end{theo}

\begin{proof}
Obviously, $\varphi(\frac{1}{a_n})\in t_\uu(\TT)$ since $u_i\varphi(\frac{1}{a_k})=0$ for $i\geq n_k$.

Let now $x\in\RR$ with $\varphi(x)\in t_\uu(\TT)$. Let $q:=\sup_{n\in\NN} q_n$, $\varepsilon:=\frac{1}{2q}$ and $n_0\in\NN$ such that 
$\Vert u_nx\Vert\leq\varepsilon$ for $n\geq n_0$. Let $k\in\NN$ with $n_k\geq n_0$. By assumption, there are $v_i\in\NN$ with $u_{n_k+i}=a_kv_i$ ($i\geq 0$). Moreover, $\frac{v_i}{v_{i-1}}=q_{n_k+i}\leq q$. Let $z=a_kx$. We show that $\Vert z\Vert=0$. Suppose that $\Vert z\Vert>0$. Then by Lemma \ref{basic1} there exists $m\in\NN$ such that $\Vert v_mz\Vert>\varepsilon$. This contradicts the fact that 
$\Vert v_mz\Vert= \Vert u_{n_k+m}x\Vert\leq\varepsilon$. We have seen that $\Vert a_kx\Vert=0$. This implies that $x\in \langle \frac{1}{a_k}\rangle$, 
equivalently $\varphi(x)\in \langle \varphi(\frac{1}{a_k})\rangle$.
\end{proof}

As mentioned in the introduction, Theorem \ref{T0} generalizes Theorems \ref{bounded} and \ref{an-dn}:
applying Theorem \ref{T0} to $\uu:=\aa$ yields \ref{bounded} and to $\uu:=\dd$ yields \ref{an-dn}; in the latter case of $\uu=\dd$ observe that $\frac{u_n}{u_{n-1}}\leq 2$ for all $n\in\NN$. 

\section{Statistically characterized subgroups of $\TT$}\label{sCh}

In this section we are interested in conditions which imply that $t_\dd^s(\TT)=t_\dd(\TT)$ where 
$\aa =(a_n)\in\A$ and $\dd =(d_n)$ is defined by ($\sharp$) of the introduction. 

Let $\aa=(a_n)\in\A$ and $q_n=\frac{a_{n}}{a_{n-1}}$ where $a_0=1$. Let $\dd=(d_n)\in\S$ be defined 
by $(\sharp)$ of the Introduction.
Then $\aa$ is a subsequence of $\dd$. We write $a_k=d_{n_k}$ for $k\in\NN$. Then $n_{k+1}=n_k+q_{k+1}-1$. 
We set $N_k:=\{n_k, n_k+1,\dots, n_{k+1}-1\}$ and $L(A):=\bigcup_{k\in A}N_k$ for $A\subseteq \NN$.
We will consider the following two conditions:

There exists a real number $\tau>1$ such that  $n_{k+1}\geq \tau n_k$ for all $k\in\NN$;
\quad\quad\quad (C1)

$\bar{d}(L(A))>0$ for any infinite $A\subseteq\NN$.\footnote{This means, in the terminology of \cite[Definition 3.9]{ADG}, that $(a_n)$ is strongly not density lifting invariant (strongly not dli for short).}
\quad\quad\quad\quad\quad\quad\quad\quad\quad\quad\quad\quad\quad\quad\quad\quad\quad (C2) 

It is easy to see that Condition (C1) implies Condition (C2) using, under the assumption (C1), the estimation 
$|N_k|/(n_{k+1}-1)\geq\frac{n_{k+1}-n_k}{n_{k+1}}= 1-\frac{n_k}{n_{k+1}} \geq 1-\frac{1}{\tau}$ 
(cf. \cite[Proposition 3.10]{ADG}).
The main result of \cite{ADG} says that Condition (C2) implies $t_\dd^s(\TT)=t_\dd(\TT)$. In \ref{T2}, we will give a very short proof of this result.
To show better the idea of the proof we first prove this in Theorem \ref{T1} under the stronger assumption (C1). This already answers the Questions \ref{Q1} and \ref{Q3} as explained in the introduction.

\begin{theo}\label{T1}
Let $\tau>1$. If $n_{k+1}\geq \tau n_k$ for all $k\in\NN$, then 
$$t_\dd^s(\TT)=t_\dd(\TT)=\varphi (\langle\{\frac{1}{a_n}:n\in\NN\}\rangle).$$ 
\end{theo}

The proof of \ref{T1} (and \ref{T2}) is based on the following lemmata.

\begin{lem}\label{L1}
Let $I\subseteq\RR$ be an interval of length $l$, $\alpha>0$ and $k:=|\alpha\ZZ\cap I|$.
Then $\frac{l}{\alpha}-1\leq k\leq\frac{l}{\alpha}+1$.
\end{lem}

\begin{proof}
This follows from $(k-1)\alpha\leq l\leq (k+1)\alpha$.
\end{proof}

\begin{lem}\label{L2} 
Let $0<\varepsilon<\frac{1}{9}$ and $0<\alpha\leq\frac{1}{2}$. Let $p\in\NN$ with $p\alpha\geq\frac{1}{4}$.

Then $|\{r\in\NN: r\leq p, \Vert r\alpha\Vert\geq\varepsilon\}|\geq\frac{p}{9}$.  
\end{lem}

\begin{proof}
Let $A=\{r\alpha: r=1,\dots,p\}$ and $\rho=|\{z\in A:\Vert z\Vert\geq \varepsilon\}|$. 

(i) Suppose first that $p\alpha\geq 1$. Let $m=[\alpha p]$. Then $m\geq 1$.

If $i\in\{0,1,\dots,m-1\}$ and $r\in\ZZ$ such that $r\alpha\in [\varepsilon,1-\varepsilon]+i$, then 
$$0<r\alpha \leq (1-\varepsilon)+(m-1)<m\leq p\alpha\,,$$ 
hence $1\leq r<p$ and $r\alpha\in A$. Therefore
$$\bigcup_{0\leq i<m} \alpha\ZZ\cap ([\varepsilon,1-\varepsilon]+i)\subseteq A\cap
\bigcup_{i\in\ZZ} \,([\varepsilon,1-\varepsilon]+i)=\{z\in A:\Vert z\Vert\geq \varepsilon\}$$
and consequently, using Lemma \ref{L1} and that $1-2\varepsilon-\alpha\geq\frac{1}{4}$, we get 
$$\rho\geq\sum_{i=0}^{m-1}|\alpha\ZZ\cap ([\varepsilon,1-\varepsilon]+i)|\geq m(\frac{1-2\varepsilon}{\alpha}-1)=p\frac{m}{p\alpha}(1-2\varepsilon-\alpha)\geq 
p\frac{m}{m+1}\frac{1}{4}\geq \frac{p}{8}$$

(ii) If $1-\varepsilon\leq p\alpha <1$, then $\{z\in A:\Vert z\Vert\geq \varepsilon\}\supseteq
\alpha\ZZ\cup [\varepsilon,1-\varepsilon]$ and therefore by Lemma \ref{L1} 
$$\rho\geq\frac{1 - 2\varepsilon}{\alpha}-1=\frac{1-2\varepsilon-\alpha}{\alpha}\geq\frac{1}{4\alpha}
=\frac{p}{4p\alpha}\geq\frac{p}{4}.$$ 

(iii) Let $\frac{1}{4}\leq p\alpha<1-\varepsilon$. If $p\leq 9$, then $\rho\geq 1\geq\frac{p}{9}$.
If $p>9$, then by Lemma \ref{L1}     
$$\rho\geq\frac{p\alpha-\varepsilon}{\alpha}-1=p\frac{(p-1)\alpha-\varepsilon}{p\alpha}\geq 
p((p-1)\alpha-\varepsilon)=p(\frac{p-1}{p}p\alpha-\varepsilon)\geq p(\frac{9}{10}\cdot\frac{1}{4}-\frac{1}{9})>\frac{p}{9}\,.$$
\end{proof}

\begin{lem}\label{L3}
Let $x\in\RR\setminus \langle\{ \frac{1}{a_i}:i\in\NN\}\rangle$, $0<\varepsilon<\frac{1}{9}$
and $N_k(\varepsilon):=\{n\in N_k: \Vert d_nx\Vert\geq\varepsilon\}$ for $k\in\NN$. 

Then for any $l\in\NN$ there exists $k\geq l$ such that $|N_k(\varepsilon)|\geq \frac{1}{9}|N_k|$.
\end{lem}

\begin{proof}
First observe that $x\in\RR\setminus \langle\{ \frac{1}{a_i}:i\in\NN\}\rangle$ implies $\Vert a_ix\Vert>0$ for all $i\in \NN$.

Let $l\in\NN$. Applying Lemma \ref{basic1} with $z=a_lx$ and $q=2$ one sees that $\Vert d_mx\Vert\geq\frac{1}{4}$ for some $m\geq n_l$. Let $k\geq l$ with $n_k\leq m<n_{k+1}$.
Then $d_m=ra_k$ with $1\leq r\leq q_{k+1}-1=:p$. Define $\alpha:=\Vert a_kx\Vert$. Then with Lemma \ref{basic0} we get 
$\frac{1}{4}\leq\Vert d_mx\Vert=\Vert r a_kx\Vert=\Vert r\alpha\Vert\leq r\alpha\leq p\alpha$.
Since $|N_k(\varepsilon)|=|\{r\in\NN: r\leq p, \Vert r\alpha\Vert\geq\varepsilon\}|$, Lemma \ref{L2} implies $|N_k(\varepsilon)|\geq\frac{p}{9}.$ Now observe that $|N_k|=n_{k+1}-n_k=q_{k+1}-1=p$.
\end{proof}

\textbf{Proof of Theorem \ref{T1}:} Obviously 
$\varphi (\langle\{\frac{1}{a_n}:n\in\NN\}\rangle)\subseteq t_\dd(\TT)  \subseteq t_\dd^s(\TT)$;
 see also Theorem \ref{an-dn}.
To show that $t_\dd^s(\TT)\subseteq\varphi (\langle\{\frac{1}{a_n}:n\in\NN\}\rangle)$
suppose that $x\in\RR\setminus \langle \{\frac{1}{a_i}:i\in\NN\}\rangle$, but $\varphi (x)\in t_\dd^s(\TT)$.

Let $\delta:=\frac{1}{10}(1-\frac{1}{\tau})$, $0<\varepsilon <\frac{1}{9}$ and $E:=\{n\in\NN: \Vert d_nx\Vert \geq \varepsilon\}$.
Then $\bar{d}(E)=0$ since $\varphi (x)\in t_\dd^s(\TT)$, and thus there exists $m_0\in\NN$ such that 
$|E\cap [1,n]|\leq n\cdot\delta$ for $n\geq m_0$. 

Define $N_k(\varepsilon)$ as in Lemma \ref{L3}. By Lemma \ref{L3} there exists $k\in\NN$  such that $n_k\geq m_0$ and  $|N_k(\varepsilon)|\geq \frac{1}{9}|N_k|$.
It follows
$$n_{k+1}\delta\geq |E\cap [1,n_{k+1}]|\geq |N_k(\varepsilon)|\geq\frac{1}{9}|N_k|=\frac{1}{9}(n_{k+1}-n_k)\,,$$
hence $\delta\geq\frac{1}{9}(1-\frac{n_k}{n_{k+1}})\geq \frac{1}{9}(1-\frac{1}{\tau})$. This contradicts the choice of $\delta$.
{~\hfill$\Box$\par \medskip} 

\begin{theo}\label{T2}(\cite[Theorem 3.15]{ADG}, \cite[Corollary 2.4]{DG24})

If $\bar{d}(L(A))>0$ for all infinite $A\subseteq\NN$, then 
$t_\dd^s(\TT)=t_\dd(\TT)=\varphi (\langle\{\frac{1}{a_n}:n\in\NN\}\rangle)$.
\end{theo}

\begin{proof}
To show that $t_\dd^s(\TT)\subseteq\varphi (\langle\{\frac{1}{a_n}:n\in\NN\}\rangle)$
suppose that $x\in\RR\setminus \langle\{ \frac{1}{a_i}:i\in\NN\}\rangle$, but $\varphi (x)\in t_\dd^s(\TT)$.

(i) Let $0<\varepsilon <\frac{1}{9}$. 
By Lemma \ref{L3} there exists $k_1<k_2<\dots$  such that $|N_{k_i}(\varepsilon)|\geq \frac{1}{9}|N_{k_i}|$ for all $i\in\NN$. Let $A:=\{k_i:i\in\NN\}$ and $0<\delta<\bar{d}(L(A))$.

(ii) Let $E:=\{n\in\NN: \Vert d_nx\Vert \geq \varepsilon\}$ and $m_0\in\NN$. We show that there exists $m> m_0$ with $|E\cap [1,n]|>\frac{1}{9}\delta m$. Let $r\in A$ with $n_r> m_0$. Since $\bar{d}(L(A))>\delta$, there exists $m>n_r$ with $|L(A)\cap [1,m]|>\delta m$. 
Let $k=\max\{l\in A :n_l\leq m\}$. One easily sees that 
$|L(A)\cap [1,m]|/m\leq |L(A)\cap [1,n_{k+1}-1]|/(n_{k+1}-1)$; therefore we may assume that $m=n_{k+1}-1$ for some $k\in A$. 
Then $E\cap [1,m]\supseteq\bigcup_{k\geq l\in A}N_l(\varepsilon)$, therefore
$$ |E\cap [1,m]|\geq\sum_{k\geq l\in A}|N_l(\varepsilon)|\geq\frac{1}{9}\sum_{k\geq l\in A}|N_l|=
\frac{1}{9}|L(A)\cap [1,m]|\geq\frac{1}{9}\delta m\,.$$

(iii) By (ii) there is a sequence $(m_i)\in\S$ such that $|E\cap [1,m_i]|\geq\frac{1}{9}\delta m_i$ for all $i\in\NN$. Therefore $\bar{d}(E)\geq\frac{1}{9}\delta$, i.e., $\bar{d}(E)>0$. 
This contradicts $\varphi(x)\in t_\dd^s(\TT)$. 
\end{proof} 

It is of interest that for the conclusion $t_\dd^s(\TT)=t_\dd(\TT)=\varphi (\langle\{\frac{1}{a_n}:n\in\NN\}\rangle)$ according to the method of proof given here the assumption (C2) seems to be very natural, exactly the same condition used in \cite{ADG}, although the proofs given here and in \cite{ADG} are completely different.

\end{document}